\documentclass[12pt,centertags,oneside]{amsart}
\usepackage{amssymb}
\usepackage{amsmath,amstext,amsthm,a4,amssymb,amscd}
\usepackage[mathscr]{eucal}
\usepackage{mathrsfs}
\usepackage{epsf}
\usepackage{hyperref}
 \usepackage{charter}
\usepackage[english]{babel}
\usepackage[T1]{fontenc}
\usepackage{typearea}

\usepackage[a4paper,width=16.2cm,top=3cm,bottom=3cm]{geometry}

\numberwithin{equation}{section}
\allowdisplaybreaks[3]

\newcommand{\field}[1]{\mathbb{#1}}

\newcommand{\R}{\field{R}}
\newcommand{\C}{\field{C}}
\newcommand{\N}{\field{N}}

 \def\cC{\mathscr{C}}
 \def\cD{\mathscr{D}}
\def\cL{\mathscr{L}}

\def\mH{\mathcal{H}}

\def\mP{\mathcal{P}}
\def\mQ{\mathcal{Q}}
\def\mR{\mathcal{R}}

\def\cF{\mathscr{F}}

\def\cL{\mathscr{L}}
\def\cP{\mathscr{P}}

 \DeclareMathOperator{\End}{End}
 \DeclareMathOperator{\Ker}{Ker}

\DeclareMathOperator{\tr}{Tr}

\newtheorem{thm}{Theorem}[section]

\theoremstyle{definition}

\theoremstyle{definition}

\newcommand{\be}{\begin{eqnarray}}
\newcommand{\ee}{\end{eqnarray}}
\newcommand{\ov}{\overline}

\newcommand{\comment}[1]{}

\begin{document}

\title{Donaldson's $Q$-operators for symplectic manifolds}
\dedicatory{In the memory of Professor Qikeng Lu}

\author{Wen Lu}

\address{School of Mathematics and Statistics,
Huazhong University of Science and Technology,
\newline\mbox{\quad}\,Wuhan, 430074, Hubei Province, China}
\email{wlu@hust.edu.cn}
\thanks{W.\ L.\ partially supported by
NNSFC 11401232}
%  Information for second author
\author{Xiaonan Ma}

\address{Institut de Math\'ematiques de Jussieu--Paris Rive Gauche,
UFR de Math{\'e}matiques,
Universit{\'e} Paris Diderot - Paris 7, Case 7012,
75205 Paris Cedex 13, France}
\email{xiaonan.ma@imj-prg.fr}
\thanks{X.\ M.\ partially supported by
NNSFC 11528103 and
funded through the Institutional Strategy of
the University of Cologne within the German Excellence Initiative}
%  Information for third author
\author{George Marinescu}

\address{Univerisit\"at zu K\"oln, Mathematisches institut,
Weyertal 86-90, 50931 K\"oln, Germany
\newline\mbox{\quad}\,Institute of Mathematics `Simion Stoilow',
Romanian Academy, Bucharest, Romania}
\email{gmarines@math.uni-koeln.de}
\thanks{G.\ M.\ partially supported by DFG funded
project SFB TRR 191}

%\date{\today}

\begin{abstract}
We prove an estimate for Donaldson's $Q$-operator
on a prequantized compact
symplectic manifold.
%previously obtained by
%Liu and Ma for K\"ahler manifolds.
This estimate is an ingredient in the recent result of Keller and Lejmi
about a symplectic generalization of Donaldson's
lower bound for the $L^2$-norm of the Hermitian scalar curvature.
   \end{abstract}

\maketitle
%\setcounter{section}{-1}
%%%%%%%%%%%%%%%%%%%%%%%%%%%
%%%%%%%%%%%%%%%%%%%%%%%%%%%%%%%%

\section{Introduction} \label{s0}
%%%%%%%%%%%%%%%%%%%%%%%%%%%%%%

The $Q$-operator
is an integral operator whose kernel is the square norm of
the Bergman kernel
of a positive line bundle (see \eqref{0.5}, \eqref{0.6}).
It was introduced by Donaldson \cite{Don09} in order
to find explicit numerical approximations
of K\"ahler-Einstein metrics on projective manifolds,
and have attracted much attention recently,
see e.g., \cite{Cao13,Joe10,Keller16,KeSe16,Liu07,Ma12}.

Using the full asymptotic expansion of the Bergman kernel
\cite{DLM06},
Liu and Ma \cite[Theorem 0.1]{Liu07} verified a statement of
Donaldson \cite[Section 4.2]{Don09} about the
relation of the asymptotics of $Q_{K_{p}}$ to the heat kernel.
Such statement was needed for the convergence of
the approximation procedure in \cite{Don09}.
In \cite{Joe10},
Liu and Ma improved the statement to
a $\cC^{m}$-estimate for
$Q_{K_{p}}$ on K\"ahler manifolds, as a crucial step towards
the result of \cite{Joe10} about the
convergence of the balancing flow to the Calabi flow.
This is a parabolic analogue of Donaldson's theorem
relating balanced embeddings to metrics
with constant scalar curvature \cite{Don01}.
Besides, such results also turn out to be important in
Cao and Keller's work \cite{Cao13} on Calabi's problem.

The purpose of the note is to extend the $\cC^{m}$-estimates
of the operators $Q_{K_{p}}$ to the
case of symplectic manifolds. This result, together with \cite{LMM17},
plays an important role in the recent work of
Keller and Lejmi \cite{Keller16} about a lower bound for
the $L^2$-norm
of the Hermitian scalar curvature. Such a lower bound was obtained
in the K\"ahler case
by Donaldson \cite{Don05}. Our proof is based on the
asymptotic expansion of the (generalized) Bergman kernel, which
in our case is the kernel of the spectral projection on lower lying
eigenstates of the normalized Bochner Laplacian.
We refer the readers
to the monograph \cite{Ma07} (see also  \cite{Ma08}, \cite{Ma10})
for more information on the Bergman kernel on symplectic manifolds.

Let us describe our result in detail. Let $(X, \omega)$
be a compact symplectic
manifold of real dimension $2n$. Let $(L, h^{L})$ be a
Hermitian line bundle on $X$, and let $\nabla^{L}$ be
a Hermitian connection on $(L, h^{L})$ with
curvature $R^{L}=(\nabla^{L})^{2}$. Let $(E, h^{E})$ be an auxiliary
Hermitian vector bundle with Hermitian connection $\nabla^{E}$.
We will assume throughout the paper that $(L, h^{L})$ satisfies
the pre-quantization condition
\begin{equation}\label{0.1}
\frac{\sqrt{-1}}{2\pi}R^{L}=\omega.
\end{equation}
We choose an almost complex structure $J$ on $TX$
(i.e., $J\in \textup{End}(TX)$ and $J^{2}=-1$) such that
$\omega$ is $J$-invariant and $\omega(\cdot, J\cdot)>0$.
The almost complex structure $J$
induces a splitting
$TX\otimes_{\R}\C=T^{(1, 0)}X\oplus T^{(0, 1)}X$,
where $T^{(1, 0)}X$ and $T^{(0, 1)}X$
are the eigenbundles of $J$ corresponding to the
eigenvalues $\sqrt{-1}$ and $-\sqrt{-1}$, respectively.

Let $g^{TX}(\cdot, \cdot):=\omega(\cdot, J\cdot)$ be
the Riemannian metric on $TX$
induced by $\omega$ and $J$.
The Riemannian volume form $dv_{X}$ of $(X, g^{TX})$ has the form
$dv_{X}=\omega^{n}/n!$.
We denote by $L^{p}:=L^{\otimes p}$ the tensor powers of $L$
for $p\in\N$
and by $h^{L^p}:=(h^L)^{\otimes p}$,
$h^{L^p\otimes E}=h^{L^p}\otimes h^E$,
the induced Hermitian metrics on $L^p$ and $L^{p}\otimes E$,
respectively.
The $L^{2}$-Hermitian product on the space
$\cC^{\infty}(X, L^{p}\otimes E)$ of smooth sections
of $L^{p}\otimes E$ on $X$ is given by
%-------------------------
\begin{equation}\label{lus0.2}
\big\langle s_{1}, s_{2}\big\rangle =
\int_{X}\big\langle s_{1}(x), s_{2}(x)\big
\rangle_{\!h^{L^p\otimes E}}\,dv_{X}(x).
\end{equation}
%-------------------------
Let $\nabla^{TX}$ be the Levi-Civita connection on $(X, g^{TX})$,
and let $\nabla^{L^{p}\otimes E}$
be the connection on $L^{p}\otimes E$ induced by $\nabla^{L}$
and $\nabla^{E}$.
Let $\{e_{k}\}$ be a local orthonormal frame of $(TX, g^{TX})$.
The Bochner Laplacian acting on $\cC^{\infty}(X, L^{p}\otimes E)$
is given by
%-------------------------
\begin{equation}
\Delta^{L^{p}\otimes E}
=-\sum_{k}\Big[\big(\nabla_{e_{k}}^{L^{p}\otimes E}\big)^{2}
-\nabla_{\nabla_{e_{k}}^{TX}e_{k}}^{L^{p}\otimes E}\Big].
\end{equation}
%-------------------------
Let $\Phi\in \cC^{\infty}(X, \textup{End}(E))$ be Hermitian
(i.e., self-adjoint with respect to $h^{E}$).
The renormalized Bochner Laplacian is defined by
\begin{align}\label{0.4a}
\Delta_{p, \Phi}=\Delta^{L^{p}\otimes E}-2\pi np+\Phi.
\end{align}
By \cite{GU}, \cite[Corollary\,1.2]{Ma02} there exists $C_{L}>0$
independent of $p$ such that
\begin{align}\label{0.4}
\textup{Spec}(\Delta_{p, \Phi})\subset [-C_{L}, C_{L}]
\cup [4\pi p-C_{L}, +\infty),
\end{align}
where $\textup{Spec}(A)$ denotes the spectrum of the operator
$A$.
Since $\Delta_{p, \Phi}$ is an elliptic operator on a
compact manifold, it has discrete spectrum and
its eigensections are smooth. Let
\begin{equation}\label{0.4a1}
\mH_{p}:=\bigoplus_{\lambda\in[-C_{L}, C_{L}]}
\Ker(\Delta_{p, \Phi}-\lambda)\subset
\cC^{\infty}(X, L^{p}\otimes E)
\end{equation}
be the direct sum of eigenspaces of $\Delta_{p, \Phi}$
corresponding to the
eigenvalues lying in $[-C_{L}, C_{L}]$.
In mathematical physics terms, the operator $\Delta_{p, \Phi}$
is a semiclassical Schr\"odinger
operator and the space $\mH_p$ is the space of its bound
states as $p\to\infty$.
 By \cite[Theorem 8.3.1]{Ma07},
\begin{align}\label{0.4b}
\textup{dim} \mH_{p}=\int_{X}\textup{Td}(T^{(1, 0)}X)
\textup{ch}(L^{p}\otimes E),
\end{align}
where $ \textup{Td}(\cdot),\textup{ch}(\cdot)$ denote
the Todd class and the Chern character
 of the corresponding complex vector bundle.
The formula (\ref{0.4b}) agrees with the Riemann-Roch-Hirzebruch
theorem and Kodaira vanishing theorem in the K\"ahler case.
The space $\mH_p$ proves to be
an appropriate replacement for the space of holomorphic
sections $H^0(X,L^p\otimes E)$ from the K\"ahler case.

Let $P_{\mH_p}$ be the orthogonal projection from
$\cC^{\infty}(X, L^p\otimes E)$ onto $\mH_p$. The kernel
$P_{\mH_p}(x, x')$ of $P_{\mH_p}$ with respect to $dv_{X}(x')$
is called a generalized Bergman kernel \cite{Ma08}.
Note that $P_{\mH_p}(x, x')\in (L^{p}\otimes E)_{x}
\otimes (L^{p}\otimes E)_{x'}^{\ast}$.
Set $\textup{Vol}(X, dv_{X})=\int_{X}dv_{X}$.
Following Donaldson \cite[\S\,4]{Don09}, we set
\begin{align}\label{0.5}
K_{p}(x, x')=\big|P_{\mH_p}(x, x')\big|^2,\ \
 R_{p}:=(\dim \mH_{p})/\textup{Vol}(X, dv_{X}).
\end{align}
Let $K_{p}, Q_{K_{p}}$ be the integral operators associated to
$K_p$ which is defined by for $f\in \cC^{\infty}(X)$,
\begin{align}\label{0.6}
(K_{p}f)(x)=\int_{X}K_{p}(x, y)f(y)dv_{X}(y),\ \
Q_{K_{p}}=\frac{1}{R_{p}}K_{p}f.
\end{align}
The operator $Q_{K_{p}}$ has been studied by
Donaldson \cite{Don09},
Liu-Ma \cite[Appendix]{Joe10}, \cite{Liu07}, and
Ma-Marinescu \cite[\S\,6]{Ma12} in the case of K\"ahler manifolds.

The main result of the note is as follows. For K\"ahler manifolds
it was obtained by Liu-Ma \cite[Appendix]{Joe10}, \cite{Liu07}.
\begin{thm}\label{t1}
For any integer $m\geqslant 0$, there exists a constant $C>0$
such that for any $f\in \cC^{\infty}(X)$,
\begin{align}\label{0.7}
\big\|Q_{K_{p}}(f)-f\big\|_{\cC^{m}(X)}\leqslant
\frac{C}{p}\big\|f\big\|_{\cC^{m+2}(X)}.
\end{align}
Moreover, \textup{(\ref{0.7})} is uniform in the following sense.
Consider $Q_{K_{p}}$ as a function of the parameters
$g^{TX}, h^{L}, \nabla^{L}, h^{E}, \nabla^{E}$ and $\Phi$,
that is,
$Q_{K_{p}}=Q_{K_{p}}(g^{TX}, h^{L}, \nabla^{L}, h^{E},
\nabla^{E},\Phi)$.
Let $\mathcal{M}$ be a subset of the infinite dimensional manifold
$\cD$ of all compatible tuples
$g^{TX}, h^{L}, \nabla^{L}, h^{E}, \nabla^{E}$ and $\Phi$.
Assume that{\,\rm:}
\begin{itemize}
\item[(i)]  the covariant derivatives in the direction $X$ of order
$\ell\leqslant 2n+m+6$
of elements of $\mathcal{M}$
form a set of tensors on $X\times\mathcal{M}$ which is bounded
in the
$\cC^0$-\,norm calculated in the direction of $\mathcal{M}${\,\rm;}
%taken with respect to the parameter $x\in X$,
\item[(ii)] the projection of $\mathcal{M}$ on the space
of Riemannian metrics is bounded below
in the $\cC^0$-\,norm.
\end{itemize}
Then there exists $C=C_{m}(\mathcal{M})$ such that \eqref{0.7}
holds for all tuples of
parameters from $\mathcal{M}$.
Moreover, the $\cC^m$-norm in \eqref{0.7} can be taken
on $X\times \mathcal{M}$.
\end{thm}

The organization of this paper is as follows. In Section \ref{s2},
we establish the asymptotic expansion of the generalized Bergman
kernel which extends \cite[\S 8.3]{Ma07}.
In Section \ref{s3}, we prove Theorem \ref{t1}.

\section{Asymptotic expansion of the generalized Bergman kernel}
\label{s2}
%%%%%%%%%%%%%%%%%%%%%%%%%%%%%%

In this section, we assume that $g^{TX}$ is an arbitrary
$J$-invariant Riemannian metric on $X$.
Let $\Delta^{L^{p}\otimes E}$ be the Bochner Laplacian acting
on $\cC^{\infty}(X, L^{p}\otimes E)$
associated with $g^{TX}$ and $\nabla^{L^{p}\otimes E}$.
Let $\Phi\in \cC^{\infty}(X, \textup{End}(E)$ be Hermitian.

Let $dv_{X}$ be the Riemannian volume form on $(X, g^{TX})$.
Now the Hermitian product on $\cC^\infty(X, L^p\otimes E)$
is induced by $h^L, h^{E}$ and $dv_{X}$.

We identify the two form $R^L$ with the Hermitian matrix
$\dot{R}^L \in \End(T^{(1,0)}X)$ such that
for $W,Y\in T^{(1,0)}X$,
\begin{equation}\label{lm4.2}
R^L (W,\ov{Y}) = \langle \dot{R}^LW, \ov{Y}\rangle.
\end{equation}
Set
\begin{align}\label{eq:2.1}
\tau= \tr|_{T^{(1,0)}X} \dot{R}^L,
\quad \mu_0=\displaystyle\inf_{{u\in T_x^{(1,0)}X,\,x\in X}}
R^L_x(u,\overline{u})/|u|^2_{g^{TX}}>0.
\end{align}
Note that if $g^{TX}= \omega(\cdot, J\cdot)$
then $\tau= 2\pi n$ and $\mu_{0}=2\pi$.

Then the renormalized Bochner Laplacian is defined as
\begin{align}\label{eq:2.4a}
\Delta_{p, \Phi}=\Delta^{L^{p}\otimes E}-\tau p+\Phi.
\end{align}
By %\cite{GU}, \cite[Corollary\,1.2]{Ma02}
the same references as in Introduction,
there exists $C_{L}>0$
independent of $p$ such that
\begin{align}\label{2.4}
\textup{Spec}(\Delta_{p, \Phi})\subset [-C_{L}, C_{L}]
\cup [2\mu_{0} p-C_{L}, +\infty),
\end{align}
Thus $\mH_{p}$ in (\ref{0.4a1}) is still well-defined
and (\ref{0.4b}) holds.

Let $P_{\mH_p}(x, x')$ be the smooth kernel of the orthogonal projection
$P_{\mH_p}$ from $\cC^\infty(X, L^p\otimes E)$
onto $\mH_p$ with respect to $dv_{X}(x')$.
In this section, we study  the asymptotics of
$P_{\mH_{p}}(x, x')$ as $p\rightarrow \infty$.

%\begin{proof}
%We study first the asymptotics of $P_{\mH_{p}}(x, x')$
%as $p\rightarrow \infty$.
Let $a^{X}$ be the injectivity radius of $(X, g^{TX})$.
We fix $\varepsilon\in (0, a^X/4)$.
Let $d(x, y)$ denote the Riemannian distance from $x$ to $y$
on $(X, g^{TX})$.
By \cite[Prop.\,8.3.5]{Ma07} and the argument
after \cite[Prop.\,8.3.5]{Ma07},
we get for any $l, m\in \mathbb{N}$ and $0<\theta<1$,
there exists $C>0$ such that
\begin{align}\label{0.27}
\Big|P_{\mH_{p}}(x, x')\Big|_{\cC^{m}(X\times X)}\leqslant Cp^{-l},
\ \
\textup{if}\ d(x, x')>\varepsilon p^{-\frac{\theta}{2}}.
\end{align}
Now we still need to understand the asymptotics of
$P_{\mH_{p}}(x, x')$ for
$d(x, x')>\varepsilon p^{-\frac{\theta}{2}}$.

We recall first the procedure of \cite[\S\,1.2]{Ma08} and
\cite[\S\,8.3]{Ma07}.

Denote by $B^{X}(x, \varepsilon)$ and
$B^{T_{x}X}(0, \varepsilon)$ the open balls
in $X$ and $T_{x}X$ with center $x$ and radius $\varepsilon$,
respectively. We identify $B^{T_{x}X}(0, \varepsilon)$ with
$B^{X}(x, \varepsilon)$
by using the exponential map of $(X, g^{TX})$.

We fix $x_{0}\in X$. For $Z\in B^{T_{x_{0}}X}(0, \varepsilon)$,
we identify $L_{Z}, E_{Z}$ and $(L^{p}\otimes E)_{Z}$
to $L_{x_{0}}, E_{x_{0}}$
and $(L^{p}\otimes E)_{x_{0}}$ by parallel transport with
respect to the connections
$\nabla^{L}, \nabla^{E}$ and $\nabla^{L^{p}\otimes E}$ along
the curve
$\gamma_{Z}: [0, 1]\ni u\rightarrow \textup{exp}^{X}_{x_{0}}(uZ)$.
Then under our identification, $P_{\mH_{p}}(Z, Z')$ is a function
on $Z, Z'\in T_{x_{0}}X$, $|Z|, |Z'|< \varepsilon$. We denote
it by $P_{\mH_{p}, x_{0}}(Z, Z')$.
Let $\pi: TX\times_{X} TX\rightarrow X$ be the natural projection
from the fiberwise
product of $TX$ on $X$. Then we can view
$P_{{\mH_{p}}, x_{0}}(Z, Z')$ as a smooth
function over $TX\times_{X} TX$ by identifying a section
$s\in \cC^{\infty}(TX\times_{X} TX, \pi^{\ast}(\textup{End}(E)))$
with the family $(s_{x})_{x\in X}$,
where $s_{x}=s|_{\pi^{-1}(x)}$.

Let $\{e_{i}\}_{i}$ be an oriented
orthonormal basis of $T_{x_{0}}X$, and let $\{e^{i}\}_{i}$ be its
dual basis. For $\varepsilon>0$ small enough, we will extend
the geometric objects from
$B^{T_{x_{0}}X}(0, \varepsilon)$ to
$\mathbb{R}^{2n}\simeq T_{x_{0}}X$
where the identification is given by
\begin{align}\label{0.8}
(Z_{1}, \ldots, Z_{2n})\in \mathbb{R}^{2n}
\longmapsto \sum_{i}Z_{i}e_{i}\in T_{x_{0}}X,
\end{align}
such that $\Delta_{p, \Phi}$ is the restriction of
a renormalized Bochner-Laplacian
on $\mathbb{R}^{2n}$ associated with a Hermitian line bundle
with positive curvature.
In this way, we replace $X$ by $\mathbb{R}^{2n}$.

At first, we denote by $L_{0}, E_{0}$ the trivial bundles with fiber
$L_{x_{0}}, E_{x_{0}}$ on
$X_{0}=\mathbb{R}^{2n}$. We still denote by
$\nabla^{L}, \nabla^{E}, h^{L}$ etc the
connections and metrics on $L_{0}, E_{0}$ on
$B^{T_{x_{0}}X}(0, 4\varepsilon)$ induced
by the above identification. Then $h^{L}, h^{E}$ is identified
to the constant
metrics $h^{L_{0}}=h^{L_{x_{0}}}, h^{E_{0}}=h^{E_{x_{0}}}$.

Let $\rho: \mathbb{R}\rightarrow [0, 1]$ be a smooth
even function such that
\begin{align}\label{0.9}
\rho(v)=1\ \ \textup{if}\ |v|<2;\ \ \rho(v)=0\ \ \textup{if}\ |v|>4.
\end{align}
Let $\varphi_{\varepsilon}: \mathbb{R}^{2n}
\rightarrow \mathbb{R}^{2n}$ is the
map defined by $\varphi_{\varepsilon}(Z)=\rho(|Z|/\varepsilon)Z$.
Then $\Phi_{0}=\Phi\circ \varphi_{\varepsilon}$ is
a smooth function on $X_{0}$.
Let $g^{TX_{0}}(Z)=g^{TX}(\varphi_{\varepsilon}(Z))$
%$J_{0}(Z)=J(\varphi_{\varepsilon}(Z))$
be the metric %and complex structure
on $X_{0}$.
Set $\nabla^{E_{0}}=\varphi_{\varepsilon}^{\ast}\nabla^{E}$.
Then $\nabla^{E_{0}}$ is the extension of $\nabla^{E}$ on
$B^{T_{x_{0}}X}(0, \varepsilon)$.
Denote by $\mR=\sum_{i}Z_{i}e_{i}=Z$
the radial vector field on $\mathbb{R}^{2n}$. We define
the Hermitian connection $\nabla^{L_{0}}$ on $(L^{0}, h^{L_{0}})$
by
\begin{align}\label{0.10}
\nabla^{L_{0}}|_{Z}=\varphi_{\varepsilon}^{\ast}\nabla^{L}
+
\frac{1}{2}(1-\rho^{2}(|Z|/\varepsilon))R_{x_{0}}^{L}(\mR, \cdot).
\end{align}
Let $R^{L_{0}}$ denote the curvature of $\nabla^{L_{0}}$
and $\{e_{i}\}_{i}$
be an orthonormal frame of $(TX_{0}, g^{TX_{0}})$.
Let $J_{0}$ be an almost complex structure on $X_{0}$ compatible with
$g^{TX_{0}}$ and $\frac{\sqrt{-1}}{2\pi}R^{L_{0}}$
such that $J_{0}=J$ on $B^{T_{x_{0}}X}(0, 2\varepsilon)$
and $J_{0}=J_{x_{0}}$ outside $B^{T_{x_{0}}X}(0, 4\varepsilon)$.
Set (cf. \eqref{eq:2.1})
\begin{align}
\tau_{0}=\frac{\sqrt{-1}}{2}\sum_{i}R^{L_{0}}(e_{i}, J_{0}e_{i}).
\end{align}
Let $\Delta^{X_{0}}_{p, \Phi_{0}}
=\Delta^{L_{0}^{p}\otimes E_{0}}-p\tau_{0}+\Phi_{0}$
be the renormalized Bochner-Laplacian on $X_{0}$ associated
to the above
data as in $(\ref{0.4a})$.
By  \cite[(1.23)]{Ma08} there exists $C_{L_{0}}>0$
such that
 \begin{align}\label{0.11}
\textup{Spec}(\Delta^{X_{0}}_{p, \Phi_{0}})\subset
[-C_{L_{0}}, C_{L_{0}}]\cup [\mu_{0} p-C_{L_{0}}, +\infty).
\end{align}
Let $S_{L}$ be an unit vector of $L_{0}$. Using $S_{L}$ and
the above
discussion, we get an isometry $L^{p}_{0}\simeq \mathbb{C}$.
Let $P_{0, \mH_p}$ be the spectral projection of
$\Delta^{X_{0}}_{p, \Phi_{0}}$
from $\cC^{\infty}(X_{0}, L^{p}_{0}\otimes E_{0})
\simeq \cC^{\infty}(X_{0}, E_{0})$ corresponding
to the interval $[-C_{L_{0}}, C_{L_{0}}]$, and let
$P_{0, \mH_p}(x, x')$ be the smooth kernel of $P_{0, \mH_p}$
with respect to
the volume form $dv_{X_{0}}(x')$.
By \cite[Proposition\,1.3]{Ma08} (for $q=0$ therein),
for any $l, m\in \mathbb{N}$, there exists $C_{l, m}>0$ such that
for $x, x'\in B^{T_{x_{0}}X}(0, \varepsilon)$, we have
\begin{align}\label{0.11a}
\Big|\big(P_{0, \mH_p}-P_{\mH_p}\big)(x, x')
\Big|_{\cC^{m}(X\times X)}\leqslant C_{l, m}p^{-l},
\end{align}
here the $\cC^{m}$-norm is induced by
$\nabla^{TX}, \nabla^{L}, \nabla^{E}, h^{L}, h^{E}$ and
$g^{TX}$.

Let $dv_{TX}$ be the Riemannian volume form on
$(T_{x_{0}}X, g^{T_{x_{0}}X})$.
Let $\kappa(Z)$ be the smooth positive function defined
by the equation
\begin{align}\label{0.12}
dv_{X_{0}}(Z)=\kappa(Z)dv_{TX}(Z),
\end{align}
with $\kappa(0)=1$. Denote by $\nabla_{U}$ the ordinary
differentiation operation
on $T_{x_{0}}X$ in the direction $U$. Denote by
$t=\frac{1}{\sqrt{p}}$.
For $s\in \cC^{\infty}(\mathbb{R}^{2n}, E_{0})$ and
$Z\in \mathbb{R}^{2n}$, set
\begin{align}\begin{split}\label{0.13}
(S_{t}s)(Z)=& s(Z/t), \ \
\nabla_{t}=tS^{-1}_{t}\kappa^{\frac{1}{2}}\nabla^{L_{0}}
\kappa^{-\frac{1}{2}}S_{t},
\\
\cL_{t}=&
S^{-1}_{t}\kappa^{\frac{1}{2}}t^{2}\Delta^{X_{0}}_{p, \Phi_{0}}
\kappa^{-\frac{1}{2}}S_{t}.
\end{split}\end{align}
It follows from (\ref{0.11}) and (\ref{0.13}) that for $t$
small enough (cf. \cite[(1.43)]{Ma08}),
\begin{align}\label{0.13a}
\textup{Spec}(\cL_{t})\subset
\big[\!-C_{L_{0}}t^{2}, C_{L_{0}}t^{2}\,\big]
\cup [\frac{1}{2}\mu_{0}, +\infty).
\end{align}
Let $\delta$ be the counterclockwise oriented circle in
$\mathbb{C}$ of center $0$ radius $\frac{1}{4}\mu_{0}$.
By (\ref{0.13a}), there exists $t_{0}>0$ such that the resolvent
$(\lambda-\cL_{t})^{-1}$ exists for $\lambda\in \delta$
and $t\in (0, t_{0}]$.

We denote by $\langle\cdot, \cdot\rangle_{0, L^2}$ and
$\|\cdot\|_{0, L^2}$ the scalar
product and the $L^2$-norm on $\cC^{\infty}(X_{0}, E_{0})$
induced by $g^{TX_{0}}$
as in (\ref{lus0.2}). For $s\in C^{\infty}(X_{0}, E_{0})$, set
\begin{align}\label{0.16c}
\|s\|^{2}_{t, 0} = &
\|s\|^{2}_{0}=\int_{\mathbb{R}^{2n}}|s(Z)|^{2}_{h^{E_{0}}}
dv_{TX}(Z),
\nonumber \\
\|s\|^{2}_{t, m} =&\sum^{m}_{l=1}\sum^{2n}_{i_{1}, \cdots, i_{l}=1}
\|\nabla_{t, e_{i_{1}}}\cdots \nabla_{t, e_{i_{l}}}s\|^{2}_{t, 0}.
\end{align}
We denote by $\langle\cdot, \cdot\rangle$ the inner product
on $C^{\infty}(X_{0}, E_{0})$
corresponding to $\|\cdot\|_{t, 0}$.
Let $H^{m}_{t}$ be the Sobolev space of order $m$ with norm
$\|\cdot\|_{t, m}$.
Let $H^{-1}_{t}$ be the Sobolev space of order $-1$ and
let $\|\cdot\|_{t, -1}$
be the norm on $H^{-1}_{t}$ defined by
$\|s\|_{t, -1}=\sup_{0\neq s'\in H^{1}_{t}}
|\langle s, s' \rangle_{t, 0}|/\|s'\|_{t, 1}$.
If $A\in \cL(H^{m}, H^{m'})$, then we denote by $\|A\|^{m, m'}_{t}$
the norm of $A$ with respect to the norms $\|\cdot\|_{t, m}$
and $\|\cdot\|_{t, m'}$.

Let $\mP_{0, t}$ the orthogonal projection from
$(\cC^{\infty}(X_{0}, E_{0}), \|\cdot\|_{0})$ onto the space of the
direct sum of eigenspaces of $\cL_{t}$ corresponding to the
eigenvalues lying in $[-C_{L_{0}}t^{2}, C_{L_{0}}t^{2}]$.
Let $\mP_{0, t}(Z, Z')=\mP_{0, t, x_{0}}(Z, Z')$ (with $Z, Z'\in X_{0}$)
be the smooth kernel of $\mP_{0, t}$ with respect to $dv_{TX}(Z')$.
We denote by $\cC^{m}(X)$ the $\cC^{m}$-norm
for the parameter $x_{0}\in X$. By \cite[(4.2.9)]{Ma07},
we have the
following extension of \cite[Theorem 1.10]{Ma08} (for $q=0$).

{\bf Claim.} For any $r, m', m\in \mathbb{N}$, there exists $C>0$
such that
for $t\in (0, t_{0}]$ and $Z, Z'\in T_{x_{0}}X$,
\begin{align}\label{0.14}
\sup_{|\alpha|+|\alpha'|\leq m'}
\bigg|\frac{\partial^{|\alpha|+|\alpha'|}}{\partial Z^{\alpha}
\partial Z'^{\alpha'}}
\frac{\partial^{r}}{\partial t^{r}}\mP_{0, t}(Z, Z')\bigg|_{\cC^{m}(X)}
\leqslant C(1+|Z|+|Z'|)^{M_{r, m', m}}
\end{align} with
\begin{align}\label{0.15}
M_{r, m', m}=2n+2+2r+m'+2m.
\end{align}

We will sketch the proof of the claim.
The readers to referred to \cite{DLM06}, \cite[Chapter 4]{Ma07}
and \cite[\S\,1]{Ma08} for more details.
In fact, by (\ref{0.13a}), for any $k\in \mathbb{N}^{\ast}$
(cf. \cite[(1.55)]{Ma08}),
\begin{align}\label{0.5a}
\mP_{0, t}=\frac{1}{2\pi \sqrt{-1}}
\int_{\delta}\lambda^{k-1}(\lambda-\cL_{t})^{-k}d\lambda.
\end{align}
For $m\in \mathbb{N}$, let $\mQ^m$ be the set of operators
$\{\nabla_{t, e_{i_1}}\cdots\nabla_{t, e_{i_j}}\}_{j\leqslant m}$.
By \cite[(1.58)]{Ma08},
\begin{align}\label{0.5b}
\big\|Q\mP_{0, t}Q'\big\|^{0, 0}_{t}\leqslant C_{m},
\ \ \textup{for}\ \ Q, Q'\in \mQ^m.
\end{align}
Let $\|\cdot\|_{m}$ be the usual Sobolev norm on
$C^{\infty}(\mathbb{R}^{n}, E_{0})$ induced by
$h^{E_{0}}$ and the volume form $dv_{TX}(Z)$.
By \cite[(4.29)]{Ma07}, there exists $C>0$ such that
for $s\in C^{\infty}(X_{0}, E_{0})$ with
$\textup{supp}(s)\subset B^{T_{x_{0}}X}(0, q)$, $m\geq 0$,
\begin{align}\label{0.5c}
\frac{1}{C}(1+q)^{-m}\|s\|_{t, m}\leqslant \|s\|_{m}
\leqslant C(1+q)^{m}\|s\|_{t, m}.
\end{align}
Now (\ref{0.5b}) and (\ref{0.5c}) together with Sobolev inequalities
imply that for $Q, Q'\in \mQ^m$,
\begin{align}\label{0.5d}
\sup_{|Z|, |Z'|\leqslant q}\Big|Q_{Z}Q'_{Z'}\mP_{0, t}(Z, Z')\Big|
\leqslant C(1+q)^{2n+2}.
\end{align}
Combining \cite[(1.35)]{Ma08} and (\ref{0.5d}) yields (\ref{0.14})
for $r=m'=0$.
To obtain (\ref{0.14}) for $r\geq 1$ and $m'=0$, note that
by (\ref{0.5a}),
\begin{align}\label{0.5e}
\frac{\partial^{r}}{\partial t^{r}}\mP_{0, t}=\frac{1}{2\pi \sqrt{-1}}
\int_{\delta}\lambda^{k-1}\frac{\partial^{r}}{\partial t^{r}}
(\lambda-\cL_{t})^{-k}d\lambda.
\end{align}
For $k, r\in \mathbb{N}^{\ast}$, let
\begin{align}\label{0.16}
I_{k, r}=\Big\{({\bf k}, {\bf r})=(k_{i}, r_{i})\big|
\sum^{j}_{i=0}k_{i}=k+j, \sum_{i=1}^{j}r_{i}=r, k_{i}+r_{i}
\in \mathbb{N}^{\ast}\Big\}.
\end{align}
Then there exist $a^{\bf k}_{\bf r}\in \mathbb{R}$ such that
\begin{align}\label{0.16b}
& A^{\bf k}_{\bf r}(\lambda, t)=(\lambda-\cL_{t})^{-k_{0}}
\frac{\partial^{r_1}\cL_{t}}{\partial t^{r_1}}
(\lambda-\cL_{t})^{-k_{1}}\cdots
\frac{\partial^{r_j}\cL_{t}}{\partial t^{r_j}}
(\lambda-\cL_{t})^{-k_{j}},
\nonumber \\
&\frac{\partial^{r}}{\partial t^{r}}(\lambda-\cL_{t})^{-k}
=\sum_{({\bf k}, {\bf r})\in I_{k, r}}a^{\bf k}_{\bf r}
A^{\bf k}_{\bf r}(\lambda, t).
\end{align}
We can now carry on nearly word by word the corresponding part
of the proof of \cite[Theorem 1.10]{Ma08}
to finish the proof of (\ref{0.14}). We finish the proof of the claim.

Set (cf. \cite[(4.1.65)]{Ma07})
\begin{align}\begin{split}\label{0.17}
\cF_{r}=&
\frac{1}{2\pi\sqrt{-1}r!}\int_{\delta}\lambda^{k-1}
\sum_{({\bf k}, {\bf r})\in I_{k, r}}a^{\bf k}_{\bf r}A^{\bf k}_{\bf r}
(\lambda, 0)d\lambda,
\\
\cF_{r, t}=&
\frac{1}{r!}\frac{\partial^{r}}{\partial t^{r}}\mP_{0, t}-\cF_{r}.
\end{split}\end{align}
Let $\cF_{r}(Z, Z')$ ($Z, Z'\in T_{x_{0}}X$) be the smooth kernel
of $\cF$ with respect to $dv_{TX}(Z')$.
Then $\cF_{r}(Z, Z')\in \cC^{\infty}(TX\times_{X} TX, \pi^{\ast}
\textup{End}(E))$. By the proof
of the estimate (\ref{0.14}), we observe that $\cF_{r}$ verifies
the similar
inequalities as (\ref{0.14}), i.e., to replace the factor
$\frac{\partial^{r}}{\partial t^{r}}\mP_{0, r}$ in (\ref{0.14})
by $\cF_{r}$.
Using this observation, (\ref{0.14}) and (\ref{0.17}),
we obtain the extension of \cite[Theorem 1.12]{Ma08}:
there exists $C>0$ such that for $t\in (0, t_{0}]$ and
$Z, Z'\in T_{x_{0}}X$,
\begin{align}\label{0.18}
\Big|\cF_{r, t}(Z, Z')\Big|\leqslant Ct^{1/(2n+1)}(1+|Z|+|Z'|)^{2n+2}.
\end{align}
By (\ref{0.17}) and (\ref{0.18}), we have (cf. \cite[(1.78)]{Ma08})
\begin{align}\label{0.19}
\frac{1}{r!}\frac{\partial^{r}}{\partial t^{r}}\mP_{0, t}|_{t=0}
=\cF_{r}.
\end{align}
By (\ref{0.14}), (\ref{0.19}) and the Taylor expansion
\begin{align}\label{0.19a}
G(t)-\sum^{k}_{r=0}\frac{1}{r!}\frac{\partial^{r}G}{\partial t^{r}}(0)
t^{r}
=\frac{1}{k!}\int^{t}_{0}(t-s)^{k}
\frac{\partial^{k+1}G}{\partial s^{k+1}}(s)ds,
\end{align}
we obtain the extension of \cite[Theorem 1.13]{Ma08}:
for any $k, m, m'\in \mathbb{N}$, there exists $C>0$ such that for
$t\in (0, t_{0}]$, $Z, Z'\in T_{x_{0}}X$ and for
$|\alpha|+|\alpha'|\leqslant m'$,
\begin{align}\label{0.20}
\bigg|
\frac{\partial^{|\alpha|+|\alpha'|}}{\partial Z^{\alpha}Z'^{\alpha'}}
\Big(\mP_{0, t}-\sum^{k}_{r=0}\cF_{r}t^{r}\Big)(Z, Z')
\bigg|_{\cC^{m}(X)}
\leqslant
Ct^{k+1}\big(1+|Z|+|Z'|\big)^{M_{k+1, m', m}}.
\end{align}
By (\ref{0.12}) and (\ref{0.13}), for $Z, Z'\in \mathbb{R}^{2n}$
(cf. \cite[(1.112)]{Ma08}),
\begin{align}\label{0.20a}
P_{0, \mH_{p}}(Z, Z')=
t^{-2n}\kappa^{-\frac{1}{2}}(Z)\mP_{0, t}(Z/t, Z'/t)
\kappa^{-\frac{1}{2}}(Z').
\end{align}
Combining (\ref{0.11a}), (\ref{0.20}) and (\ref{0.20a}), we obtain
\begin{align}\begin{split}\label{0.21}
 \bigg|&
\frac{\partial^{|\alpha|+|\alpha'|}}{\partial Z^{\alpha}Z'^{\alpha'}}
\Big(\frac{1}{p^{n}}P_{\mH_p, x_{0}}(Z, Z')
-\sum^{k}_{r=0}\cF_{r}(\sqrt{p}Z, \sqrt{p}Z')
\kappa^{-\frac{1}{2}}(Z)\kappa^{-\frac{1}{2}}(Z')
p^{-\frac{r}{2}}\Big)\bigg|_{\cC^{m}(X)}\\ &
\leqslant Cp^{-\frac{k-m+1}{2}}
\big(1+\sqrt{p}|Z|+\sqrt{p}|Z'|\big)^{M_{k+1, m', m}}.
\end{split}\end{align}
Now we fix $k_{0}, m', m$. Take
\begin{align}\label{0.21a}
 k=k_{0}+m'+2\ \ \textup{and}\ \
\theta=1/\big(2M_{k+1, m', m}\big).
\end{align}
Then for $|\alpha|+|\alpha'|\leqslant m'$
and $|Z|, |Z'|< p^{-\frac{1}{2}+\theta}$, we have
\begin{align}\begin{split}\label{0.22}
 \bigg|&
\frac{\partial^{|\alpha|+|\alpha'|}}{\partial Z^{\alpha}Z'^{\alpha'}}
\Big(\frac{1}{p^{n}}P_{\mH_p, x_{0}}(Z, Z')
-\sum^{k}_{r=0}\cF_{r}(\sqrt{p}Z, \sqrt{p}Z')
\kappa^{-\frac{1}{2}}(Z)
\kappa^{-\frac{1}{2}}(Z')p^{-\frac{r}{2}}\Big)
\bigg|_{\cC^{m}(X)}
\\
 &\leqslant
Cp^{-\frac{k_{0}}{2}-1}.
\end{split}\end{align}
To sum up, we have finished the proof of the following result:

\begin{thm} \label{t2.1}
 For
any $k_{0}, m', m\in \mathbb{N}$, there exists $C>0$ such that
for $|\alpha|+|\alpha'|\leqslant m'$ and
$|Z|, |Z'|<p^{-\frac{1}{2}+\theta}$ with
\begin{align}\label{0.22a}
\theta=\frac{1}{2\big(2n+8+2k_{0}+3m'+2m\big)},
\end{align}
we have
\begin{align}\begin{split}\label{0.23}
& \bigg|
\frac{\partial^{|\alpha|+|\alpha'|}}{\partial Z^{\alpha}Z'^{\alpha'}}
\Big(\frac{1}{p^{n}}P_{\mH_p, x_{0}}(Z, Z')
-\sum^{k}_{r=0}\cF_{r}(\sqrt{p}Z, \sqrt{p}Z')
\kappa^{-\frac{1}{2}}(Z)
\kappa^{-\frac{1}{2}}(Z')p^{-\frac{r}{2}}\Big)
\bigg|_{\cC^{m}(X)} \leqslant
Cp^{-\frac{k_{0}}{2}-1},
\end{split}\end{align}
where $k=k_{0}+m'+2$.
\end{thm}

We choose $\{w_{j}\}_{j=1}^n$ an orthonormal
basis of $T_{x_{0}}^{(1, 0)}X$ such that
\begin{align}\label{eq:2.26}
\dot{R}^L_{x_{0}}={\rm diag}
(a_{1},\cdots, a_{n})\in \End( T_{x_{0}}^{(1, 0)}X).
\end{align}
Then
$e_{2j-1}=\tfrac{1}{\sqrt{2}}(w_j+\overline{w}_j)$ and
$e_{2j}=\tfrac{\sqrt{-1}}{\sqrt{2}}(w_j-\overline{w}_j)\,,
 j=1,\dotsc,n$,
form an orthonormal basis of $T_{x_0}X$.
We use the coordinates on $T_{x_0}X\simeq\R^{2n}$ induced by
$\{ e_i\}$ as in \eqref{0.8}
and in what follows we also introduce
the complex coordinates $z=(z_1,\cdots,z_n)$
on $\C^n\simeq\R^{2n}$.
Set
%(cf. \cite[(1.91)]{Ma08} with $a_{j}=2\pi$ therein)
\begin{align}\label{0.24}
\cP(Z, Z')=\textup{exp}\Big[-\frac{1}{4}
\sum^{n}_{j=1}a_{j}\big(|z_{j}|^{2}+|z'_{j}|^{2}
- 2z_{j}\overline{z}'_{j} \big)\Big].
\end{align}

By \cite[Theorem 1.18]{Ma08}, there exist $J_{r}(Z, Z')$ polynomials
in $Z, Z'$
with the same parity as $r$  and degree $\leq 3r$ such that
\begin{align}\label{0.25}
\cF_{r}(Z, Z')=J_{r}(Z, Z')\cP(Z, Z'),\ \ J_{0}(Z, Z')=1.
\end{align}

\section{Proof of Theorem \ref{t1}} \label{s3}
%%%%%%%%%%%%%%%%%%%%%%%%%%%%%%

Now $g^{TX}(\cdot, \cdot):=\omega(\cdot, J\cdot)$,
thus $a_{j}=2\pi$ in (\ref{0.24}).

Recall that the classical heat kernel on $\mathbb{C}^{n}$ is
$e^{-u\Delta}(Z, Z')=(4\pi u)^{-n}e^{-\frac{1}{4u}|Z-Z'|^{2}}$. Then
\begin{align}\label{0.24a}
\big|\cP(Z, Z')\big|^{2}=e^{-\pi|Z-Z'|^{2}}
=e^{-\frac{\Delta}{4\pi}}(Z, Z').
\end{align}
Note that
$\big|P_{\mH_{p}, x_{0}}(Z, Z')\big|^{2}=
P_{\mH_{p}, x_{0}}(Z, Z')\overline{P_{\mH_{p}, x_{0}}(Z, Z')}$.
By (\ref{0.5}), (\ref{0.23}),  (\ref{0.25}) and (\ref{0.24a}),
there exist polynomials $J'_{r}(Z, Z')$ in $Z, Z'$ such that
for $|Z|,|Z'|< p^{-\frac{1}{2}+\theta}$
with $\theta$ in (\ref{0.22a}),
\begin{align}\label{0.26}
 \bigg|
\frac{1}{p^{2n}}K_{p, x_{0}}(Z, Z')
-\Big(1+\sum^{k}_{r=1}p^{-\frac{r}{2}}
J'_{r}(\sqrt{p}Z, \sqrt{p}Z')\Big)e^{-\pi p|Z-Z'|^{2}}
\bigg|_{\cC^{m}(X)}
\leqslant
Cp^{-\frac{k_{0}}{2}-1},
\end{align}
with
\begin{align}\label{0.26a}
J'_{1}(0, Z')=(J_{1}+\overline{J_{1}})(0, Z').
\end{align}

For a function $f\in \cC^{\infty}(X)$, we denote by
$f_{x_{0}}(Z)$ the function $f$ in normal coordinates $Z$ around
a point
$x_{0}\in X$. We have thus a family $(f_{x_{0}})$ of functions
indexed by the parameter $x_0\in X$.
Combining (\ref{0.5}), (\ref{0.27})
with $\theta$ in (\ref{0.22a}), and (\ref{0.26}),  we obtain
\begin{align}\begin{split}\label{0.28}
 \bigg|&
\frac{1}{p^{n}}K_{p}f
-p^{n}\int_{|Z'|\leqslant \varepsilon p^{-\theta/2}}
\Big(1+\sum^{k}_{r=1}p^{-\frac{r}{2}}J'_{r}(0, \sqrt{p}Z')\Big)
e^{-\pi p|Z'|^{2}}f_{x_{0}}(Z')dv_{X}(Z')
\bigg|_{\cC^{m}(X)}
\\ &\leqslant
Cp^{-\frac{k_{0}}{2}-1}\big|f\big|_{\cC^{m}(X)}.
\end{split}\end{align}

% Let $\Delta$ be the positive Laplace operator on $(X, g^{TX})$
% acting on functions on $X$,
% and let $e^{-t\Delta}(x, x')$ be the smooth kernel of
% the heat operator $e^{-t\Delta}$ with
% respect to $dv_{X}(x')$. By \cite[(35)]{Liu07}, there exist
% $\phi_{i, x_{0}}(Z')$ such
% that uniformly for $x_{0}\in X$, $Z'\in T_{x_{0}}X$ with
% $|Z'|\leqslant \varepsilon$, we have
% the following asymptotic expansion when $u\rightarrow 0$:
% \begin{align}\label{1}
% \Big|e^{-u\Delta}(0, Z')-(4\pi u)^{-n}
% \big(1+\sum^{k}_{i=1}u^{i}\phi_{i, x_{0}}(Z')\big)
% e^{-\frac{1}{4u}|Z'|^{2}}\Big|_{\cC^{0}(X)}
% =O(u^{k-n+1}).
% \end{align}
% By \cite[(2.8)]{BGV}, when $u\rightarrow 0$,
% \begin{align}\label{2}
% \Big|e^{-u\Delta}f-\sum^{k}_{i=1}\frac{(-u\Delta)^{i}}{i!} f
% \Big|_{\cC^{m}(X)}=O(u^{k+1}).
% \end{align}
% By (\ref{1}) (with $u=\frac{1}{4\pi p}$) and (\ref{2}),

By using Taylor expansion of $f_{x_{0}}(Z')$ at $0$, we obtain
\begin{align}\begin{split}\label{0.29}
& \bigg|p^{n}\int_{|Z'|\leqslant \varepsilon p^{-\theta/2}}
J'_{r}(0, \sqrt{p}Z')e^{-\pi p|Z'|^{2}}f_{x_{0}}(Z')dv_{X}(Z')
\bigg|_{\cC^{m}(X)}
 \leqslant
C\big|f\big|_{\cC^{m}(X)},
\\ &
\bigg|p^{n}\int_{|Z'|\leqslant \varepsilon p^{-\theta/2}}
e^{-\pi p|Z'|^{2}}f_{x_{0}}(Z')dv_{X}(Z')-f(x_{0})
\bigg|_{\cC^{m}(X)}
 \leqslant
\frac{C}{p}\big|f\big|_{\cC^{m+2}(X)}.
\end{split}\end{align}
Finally, by \cite[Theorem 1.18]{Ma08} and
\cite[(1.97),\,(1.98),\,(1.111)]{Ma08}, we obtain
\begin{align}\begin{split}\label{0.29a}
\int_{Z'\in \mathbb{C}^{n}}&\overline{J_{1}}(0, Z')|\cP|^{2}(0, Z')dZ'
\\=&
\int_{Z'\in \mathbb{C}^{n}}\cP(0, Z')J_{1}(Z', 0)\cP(Z', 0)dZ'
\\=&
(\cP J_{1}\cP)(0, 0)=0.
\end{split}\end{align}
Combining %(\ref{1}), (\ref{2})
Taylor expansion of $f_{x_{0}}(Z')$ at $0$, and (\ref{0.29a}) yields
\begin{align}\label{0.29b}
\bigg|p^{n}\int_{|Z'|\leqslant \varepsilon p^{-\theta/2}}
p^{-1/2}J'_{1}(0, \sqrt{p}Z')e^{-\pi p|Z'|^{2}}f_{x_{0}}(Z')dv_{X}(Z')
\bigg|_{\cC^{m}(X)}
 \leqslant
\frac{C}{p}\big|f\big|_{\cC^{m+2}(X)}.
\end{align}
Combining (\ref{0.28}) for $k_{0}=0$, (\ref{0.29})
and (\ref{0.29b}) yields
\begin{align}\label{0.30}
 \Big|
\frac{1}{p^{n}}K_{p}f-f
\Big|_{\cC^{m}(X)}
 \leqslant  \frac{C}{p}\big|f\big|_{\cC^{m+2}(X)}.
\end{align}
Then the desired $\cC^{m}$-estimate (\ref{0.7}) follows from
(\ref{0.6}) and (\ref{0.30}).
The proof of the uniformity assertion from Theorem \ref{t1} is
modeled on
\cite[\S\,4.1.7]{Ma07}, \cite[\S\,1.5]{Ma08}.
First we notice that in the proof of the estimate
(\ref{0.14}), we only use the derivatives of the coefficients of
$\cL_{t}$ with order $\leqslant 2n+m+m'+r+2$.
Thus, by (\ref{0.19a}), the constants in (\ref{0.14}), (\ref{0.18})
(resp. (\ref{0.20}), (\ref{0.21}))
are bounded, if with respect to a fixed metric $g_{0}^{TX}$,
the $\cC^{2n+m+m'+r+3}$ (resp. $\cC^{2n+m+m'+k+4}$)-norms
on $X$ of the data
$g^{TX}, h^{L}, \nabla^{L}, h^{E}, \nabla^{E}$ and $\Phi$
are bounded and $g^{TX}$ is bounded below.
Note $k=k_{0}+m'+2$ in (\ref{0.23}).
Then the constants in (\ref{0.23}) (resp. (\ref{0.26}), (\ref{0.28}),
(\ref{0.30}))
are bounded if with respect to a fixed metric $g_{0}^{TX}$,
the $\cC^{2n+m+2m'+k_{0}+6}$ (resp. $\cC^{2n+m+k_{0}+6}$,
$\cC^{2n+m+k_{0}+6}$, $\cC^{2n+m+6}$)-norm on $X$ of the data
$g^{TX}, h^{L}, \nabla^{L}, h^{E}, \nabla^{E}$ and $\Phi$ are
bounded and $g^{TX}$ is bounded below.
Moreover, taking derivatives with respect to the
parameters we obtain a similar equation to
(\ref{0.5e}) (cf. \cite[(1.65)]{Ma08}). Thus
the $\cC^{m}$-norm in (\ref{0.30}) can also include the parameters
of the $\cC^{m}$-norm
if the $\cC^{m}$-norms (with respect to the
parameter $x_{0}\in X$) of derivatives of the above data with
order $\leqslant 2n+6$ are bounded.
Thus we can take $C$ in (\ref{0.7}) independent of $g^{TX}$.
The proof of Theorem \ref{t1} is completed.
%\end{proof}

\end{document}